\documentclass[11pt]{article}
\usepackage{amsfonts}
\usepackage{mathrsfs}
\usepackage{amsmath}
\usepackage{amssymb}
\usepackage{graphicx}
\graphicspath{{figure/}}
\usepackage{epic}
\usepackage{amsthm,latexsym}
\renewcommand{\paragraph}{\roman{paragraph}}
\newtheorem{theorem}{\scshape \mdseries  Theorem}[section]
\newtheorem{lemma}[theorem]{\scshape \mdseries  Lemma}

\newtheorem{prop}[theorem]{\scshape \mdseries  Proposition}

\begin{document}

\title{\sf Finite groups admitting  a regular tournament\\ $m$-semiregular representation}
\author{Dein Wong\thanks{Corresponding author, E-mail address:wongdein@163.com.    Supported by the National Natural Science Foundation of China (No.12371025).},\ \ \  Songnian Xu, \ \ Chi Zhang 
 \\ {\small  \it School of Mathematics, China University of Mining and Technology, Xuzhou,  China.} \\Jinxing Zhao \thanks{E-mail address:zhjxing@imu.edu.cn. 
  Supported by the National Natural Science Foundation of China (No.12161062)}
\\ {\small  \it  School of Mathematics Sciences, Inner Mongolia University, Hohhot, China.}\\
 }

\date{}
\maketitle

\noindent {\bf Abstract:}\ \ For a positive integer $m$, a finite group $G$ is said to admit a
tournament $m$-semiregular representation (TmSR for short) if there exists a tournament $\Gamma$ such that the automorphism
group of $\Gamma$ is isomorphic to $G$ and acts semiregularly on the
vertex set of $\Gamma$ with $m$ orbits. Clearly, every finite group of even order does not admit a TmSR for any
positive integer $m$, and T1SR is the well-known tournament
regular representation (TRR for short).  In 1986, Godsil \cite{god} proved, by a probabilistic approach, that the only finite groups of odd order
without a TRR are $\mathbb{Z}_3^2$ and $\mathbb{Z}_3^3$ .
  More recently, Du \cite{du} proved that every finite group of
odd order has a TmSR for every $m \geq 2$. The author of \cite{du} observed that a finite group of odd order has no regular TmSR when $m$ is an even integer,   a group of order $1$ has no regular T3SR,   and $\mathbb{Z}_3^2$ admits a regular T3SR.  At the end of \cite{du}, Du proposed the following problem.

\noindent{\sf\it Problem.} \ \ {\it For every odd integer $m\geq 3$, classify finite groups of odd order which have
a regular TmSR.}

The motivation of this paper is to give an answer for the above problem. We proved that if $G$ is a finite group with odd order $n>1$, then $G$ admits a regular TmSR for any odd  integer $m\geq 3$.


 \section{Introduction}
\quad \quad
Let $G$ be a permutation group on a set $X$. For $x\in X$, denote by $G_x$ the stabilizer
of $x$ in $G$. The group $G$ is said to be semiregular on $X$ if
$G_x= \{e\}$ for every $x\in X$, and regular if it is semiregular and transitive.

A classic problem asked by K$\ddot{o}$nig \cite{kon} in 1936 is: Whether a given group can be represented as the automorphism group of a graph. K$\ddot{o}$nig's publication \cite{kon} is considered as the start of graph
representation. The above problem was solved by Frucht \cite{fru} in 1939, and he improved the answer by showing that every finite group is the automorphism group of a cubic graph in
\cite{fru}. A group $G$ is said to admit a digraphical or a graphical regular representation (DRR or
GRR for short) if there exists a digraph or a graph $\Gamma$ such that  Aut$(\Gamma)\cong G$ is regular on
the vertex set, respectively.  The most natural question is the ``GRR and DRR
problem": {\it Which group admits a GRR  or a DRR?}
Babai (\cite{bab1}, Theorem 2.1) proved that, except for
$ Q_8, \mathbb{Z}_2^2, \mathbb{Z}_2^3, \mathbb{Z}_2^4, \mathbb{Z}_3^2$,
every group admits a DRR.
Spiga \cite{spi} classified the finite groups admitting an oriented regular representation.
It is clear that if a group admits a GRR then it also admits a DRR, but the converse is
not true. After a long series
of partial results by various authors, the classification of groups with a GRR was completed
by Godsil in \cite{god1}.
A group $G$ is said to admit  a    GmSR (resp.,  DmSR) if there exists a   graph (resp., a digraph) $\Gamma$ such
that ${\rm Aut}(\Gamma) \cong G$  is semiregular on   $V(\Gamma)$ with $m$ orbits. Finite groups   admiting a GmSR or a DmSR are classified by Du et al. in \cite{du1}.


 A tournament is a digraph $\Gamma = (V, A)$  with vertex set $V$ and arc set $A$
such that, for every two distinct vertices $x,y\in V$, exactly one of $(x, y)$ and $(y, x)$ is in $A$.
It is clear that the tournament has no automorphism of order $2$.  For a (finite or infinite) group without elements of order
2, Moon \cite{moo} and Babai \cite{bab} proved that there exists a tournament whose automorphism
group is isomorphic to $G$. A group $G$ has a tournament regular representation (TRR for
short) if there exists a tournament $\Gamma$ such that Aut$(\Gamma)\cong G$ is regular on the vertex set of $\Gamma$. Babai and Imrich \cite{bab2} proved that every group of odd order except for $\mathbb{Z}_3^2$
admits a tournament regular representation. Godsil \cite{god} corrected the result by showing that  $\mathbb{Z}_3^2$ and $\mathbb{Z}_3^3$
are the only finite groups of odd order without   TRRs. We refer the reader to \cite{du2}, \cite{imr}, [13-20] for more information on graphical   representation of finite groups.


 A group $G$ is said to admit a tournament $m$-semiregular representation (TmSR for
short) for a positive integer $m$, if there exists a tournament $\Gamma$ such that  Aut($\Gamma)\cong   G$
is semiregular on the vertex set of $\Gamma$  with $m$ orbits. Note that T1SR is exactly TRR,
and recall that every finite group of odd order, except for  $\mathbb{Z}_3^2$ and $\mathbb{Z}_3^3$,
 has a TRR.  For the case when $m=2$, Babai and Imrich \cite{bab2} proved that if a group $G$ without elements of
order 2, and if $ G$ has an independent generating set, then $G$ has a T2SR. Recently, Du \cite{du} proved the following theorem:
\begin{prop} {\rm (\cite{du}, Theorem 1.2)} \ Let $m\geq 2$ be an integer. Then a finite group has a TmSR if and only if it
has odd order.\end{prop} The author of \cite{du}  observed that a finite group of odd order has no regular TmSR when $m$ is an even integer,   a group of order $1$ has no regular T3SR,  and $\mathbb{Z}_3^2$ admits a regular T3SR.  Du terminated \cite{du} with the following problem:

\noindent{\sf\it Problem.} \ \ {\it For every odd integer $m\geq 3$, classify finite groups of odd order which has
a regular TmSR.}

The motivation of this paper is to give an answer for the above problem. We obtain a result as follows.
\begin{theorem}  Let $G$ be a finite group of order $n>1$. Then $G$ admits a regular TmSR, with $m>1$, if and only if $n$ and $m$ are   odd  integers.\end{theorem}

Du  \cite{du} has indicated that, if $n$ or $m$ is even, then $G$ fails to admit a regular TmSR. Thus, we only need to prove the sufficiency for the theorem.

\section{Preliminaries and notations}
\quad\quad For a digraph $D$ with  a vertex $u$, let $N_D^+(u)=\{v\in D \ | \ (u,v)\in A(D)\}$, called the set of out-neighbors of $u$ in $D$. Set $d_D^+(u)=|N_D^+(u)|$, called the out-degree of $u$ in $D$. The in-degree of $u$ in $D$ is defined as $d_D^-(u)=|N_D^-(u)|$, where  $N_D^-(u)=\{v\in D \ |\  (v,u)\in A(D)\}$ is the set of in-neighbors of $u$ in $D$.

 We recall two known results from \cite{god} and \cite{du} for later use.
\vskip 2mm
\begin{prop} {\rm (\cite{god}, Corollary 4.4)}\ \ If $G$ is not isomorphic to $\mathbb{Z}_3^2$ or $\mathbb{Z}_3^3$, then $G$ has a TRR.
 \end{prop}
 This proposition tell us that:  If $G$ is not isomorphic to $\mathbb{Z}_3^2$ or $\mathbb{Z}_3^3$, then there is a subset $R$ of $G$  with $R\cap R^{-1}= \emptyset$ and  $R \cup R^{-1}=G\setminus\{e\}$,
   such that  Aut$({\rm Cay}(G, R))\cong G$ is  regular on $G$.
 \begin{prop} {\rm (\cite{du}, Page 5)}\ \ Suppose   $G = \langle a, x\rangle \cong  \mathbb{Z}_3^2$.   Let $R = \{a, x, ax, a^2 x\}$ and $S = \{a, x^2 , ax, a^2 x\}$ be two
subsets of $G$, where $R \cap R^{-1} = S \cap S^{-1} =\emptyset $, $R \cup R^{-1}=S \cup S^{-1}=G\setminus\{e\}$. Then $G$ has a T2SR of the form ${\rm Cay}(G, T_{ij} :\ i,j\in \mathbb{Z}_2)$, where $T_{1,1}=R$, $T_{2,2}=S$, $ T_{1,2}=\{e\}$ and $T_{2,1}=G\setminus\{e\}$.
\end{prop}

\begin{prop}  {\rm (\cite{du}, Page 5)}\ \ If $G = \langle a, x, c\rangle \cong  \mathbb{Z}_3^3$. Let $$R = \{a,x, c, ax, ac,xc,axc,a^2 x,a^2c,x^2 c,a^2 xc,a^2 x^2 c,
ax^2 c\},$$   $$S = \{a,x^2, c^2, ax, ac,xc,axc,a^2 x,a^2c,x^2 c,a^2 xc,a^2 x^2 c,
ax^2 c\}$$  be two subsets of
$G$, where $R \cap R^{-1} = S \cap S^{-1} =\emptyset $,  $R \cup R^{-1} = S \cup S^{-1} =G\setminus\{e\}$. Then $G$ has a T2SR of the form ${\rm Cay}(G, T_{i,j} :\ i,j\in \mathbb{Z}_2)$, where $T_{1,1}=R$, $T_{2,2}=S$, $ T_{1,2}=\{e\}$ and $T_{2,1}=G\setminus\{e\}$.
\end{prop}

Indeed, Du proved an extensive result: If  $G \cong  \mathbb{Z}_3^2$ or $G \cong  \mathbb{Z}_3^3$, then $G$ has a TmSR  for an arbitrary integer $m\geq 2$. In the present article,  we  only apply the assertion that ``$G$ has a T2SR" to prove our main result.

Suppose $G$ is a finite group with odd order $n>1$ and $m=2k+1\geq 3$. The identity element of $G$ is written as $e$. In order to prove that $G$ admits a regular TmSR, we  construct a  regular $m$-Cayley digraph $\Gamma={\rm Cay}(G, T_{i,j} : \ i, j\in \mathbb{Z}_m)$ such that $\Gamma$ is a tournament,
then we try to prove   that ${\rm Aut}(\Gamma)\cong G$ is semiregular on $V(\Gamma)$ with $m$ orbits. The vertex set and arc set of the above  $m$-Cayley digraph $\Gamma$ are defined as:
  $$V(\Gamma)=\bigcup_{i\in \mathbb{Z}_m} G_i,\ \ A(\Gamma)=\bigcup_{i,j\in \mathbb{Z}_m} \{(g_i, (tg)_j)\  |\ t\in T_{i,j}\},$$ where $G_i=\{g_i\ |\ g\in G\}$ and $T_{i,j}$  are carefully selected  subsets of $G$ in different way according to different   $m$ and different cases of $G$.

 Suppose  that $\Gamma={\rm Cay}(G, T_{i,j} : \ i, j\in \mathbb{Z}_m)$ is a regular tournament constructed as above. Any given element $g\in G$ induces an automorphism $\rho_g$ of $\Gamma$ in the way $f_i\mapsto (fg)_i, \forall \ f\in G, i\in \mathbb{Z}_m$, which is called the right regular representation with respect to $g$.
 All such automorphisms $\rho_g$, with $g\in G$, form  a subgroup of  ${\rm Aut}(\Gamma)$, which is denoted by ${\rm Reg}(G)$. Indeed, the action of ${\rm Reg}(G)$ on $V(\Gamma)$ is   semiregular   with $m$ orbits and $G$ is isomorphic to ${\rm Reg}(G)$. To  prove that $G$ has a regular  TmSR, it suffices to prove that ${\rm Aut}(\Gamma)={\rm Reg}(G)$ for some $\Gamma$.
When we concentrate on proving ${\rm Aut}(\Gamma)={\rm Reg}(G)$, one key technique is to show that   ${\rm Aut}(\Gamma)$ has $m$ orbits on $V(\Gamma)$. Particularly, we need to distinguish one $G_i$ from other $G_j$ with $j\not=i$ (if $G$ is   isomorphic to $\mathbb{Z}_3^2$ or $\mathbb{Z}_3^3$, we try to prove that   ${\rm Aut}(\Gamma)$ stabilizes some  $G_i\cup G_{i+k}$).

 The proof of Theorem 1.2 will be partitioned into three parts corresponding to three cases: Case I: $m\geq 7$ or $m=5, n\geq 5$.  Case II:  $m=5$ and $n=3$. Case III: $m=3$.

\section{Case I:    $m\geq 7$,  or $m=5$ and $ n\geq 5$.}
\quad\quad If  $m\geq 7$ or $m=5, n\geq 5$, $T_{i,j}$ are chosen in the following way for the definition of $\Gamma$.
\vskip 2mm
\noindent $\bullet$\ \ If $G$ admits a TRR,  we set  $T_{i,i}=R$ for  all  $i\in  \mathbb{Z}_m$ (where $R$ is given in Proposition 2.1) and  for a given $i\in \mathbb{Z}_m,$  except for $T_{0,1}$ and $T_{1,0}$,  we set $T_{i,i+j}=\{e\}$,  $T_{ i+j,i}=G\setminus\{e\}$ for $ j=1, \ldots, k$; and $T_{i,i+j}=G\setminus\{e\}$, $T_{ i+j,i}= \{e\}$ for $j=k+1, \ldots, 2k$;   particularly, we set $T_{0,1}=\{x\}, T_{1,0}=G\setminus\{x^{-1}\}$, where $x$ is any given element in $R$
\vskip 2mm
 \noindent $\bullet$\ \
  If   $G = \langle a, x\rangle \cong  \mathbb{Z}_3^2$ or $G= \langle a, x, c\rangle \cong  \mathbb{Z}_3^3$,  we set $T_{k+1,k+1}=S$ and $T_{i,i}=R$ for each other  $i\in  \mathbb{Z}_m$ with $i\not=k+1$ (where $R$ and $S$ are given in Proposition 2.2 and 2.3);  for a given $i\in \mathbb{Z}_m,$  except for $T_{0,1}$ $T_{1,0}$,  we set $T_{i,i+j}=\{e\}$, $T_{ i+j,i}=G\setminus\{e\}$ for $ j=1, \ldots, k$; and $T_{i,i+j}=G\setminus\{e\}$, $T_{ i+j,i}= \{e\}$ for $j=k+1, \ldots, 2k$;    particularly, we set $T_{0,1}=\{x\}$,
    $T_{1,0}=G\setminus\{x^{2}\}$.
 \vskip 2mm

 The graph $\Gamma$ is as just defined.  For $i\in \mathbb{Z}_m$, let $\Gamma_i$ be the sub-digraph of $\Gamma$ induced by $N_\Gamma^+(e_i)$. Denote $G_j\cap N_{\Gamma}^+(e_i)$ by $W_{i,j}$ with $j\in\mathbb{ Z}_m$.  The vertex set of $\Gamma_i$ has a partition as: $$V(\Gamma_i)=\bigcup_{j\in\mathbb{ Z}_m} W_{i,j}.$$  The vertices of $W_{i,j}$ are as follows.
\vskip 2mm
\noindent $\bullet$\ \ $W_{i,i}=S_i$ if $G\cong\mathbb{Z}_3^2$ or $\mathbb{Z}_3^3$ and $i=k+1$;    $W_{i,i}=R_i$ for all other cases.

\noindent $\bullet$\ \ $W_{i,i+1}=\{z_{i+1}\}$, where $z=x$ if $i=0$ and $z=e$ if $i\not=0$;

 \noindent $\bullet$\ \ $W_{i,j}=\{e_j\}$ for $i+2\leq j\leq i+k$;

\noindent $\bullet$\ \ $W_{i,s} =(G\setminus\{e\})_s$ for $i+k+1\leq  s\leq  i+2k-1$;

\noindent $\bullet$\ \ If   $i\not=1$, then $W_{i,i+2k}=(G\setminus\{e\})_{i+2k}$;   if $i=1$  then $W_{i,i+2k}=(G\setminus\{x^{-1}\})_{i+2k}$.
\vskip 2mm
  When $G$ admits a TRR, we  try to distinguish $G_1$  from other $G_j, j\not=1$. To achieve the goal    we try to prove that the in-degrees of vertices  in $\Gamma_1$ do not coincide with the in-degrees of vertices  in other $\Gamma_i$. However, we do not     consider the vertices in $W_{i,i}$, since the set   $\Delta_i=\{d_{\Gamma_i}^-(v)|\ v\in W_{i,i}\}$ does not depend on $i$ (see Lemma 3.1).
 \begin{lemma} For $i\not=l\in \mathbb{Z}_m$, if $T_{i,i}=T_{l,l}$, then $\Delta_i=\Delta_l$.
\end{lemma}
\begin{proof} The assumption $T_{i,i}=T_{l,l}$  implies that  $W_{i,i}=R_i $ and $W_{l,l}=R_l$.  Suppose $w_i\in W_{i,i}$ with $w\in R$. It is easy to see that, for each $i+1\leq j\leq i+2k$, exactly one vertex in $W_{i,j}$ is an in-neighbor of $w_i$ in $\Gamma_i$.
Thus, $ d_{\Gamma_i}^-(w_i)=2k+|N_{\Gamma_i}^-(w_i)\cap W_{i,i}|$. Similarly, $ d_{\Gamma_l}^-(w_l)=2k+|N_{\Gamma_l}^-(w_l)\cap W_{l,l}|$. Since $W_{i,i}=R_i $ and $W_{l,l}=R_l$, it is easy to see that $|N_{\Gamma_i}^-(w_i)\cap W_{i,i}| =|N_{\Gamma_l}^-(w_l)\cap W_{l,l}|$. Consequently,  $ d_{\Gamma_i}^-(w_i)=d_{\Gamma_l}^-(w_l)\in \Delta_l$ and $\Delta_i\subseteq \Delta_l$. Similarly,  $\Delta_l\subseteq \Delta_i$, forcing  $\Delta_i=\Delta_l$.
\end{proof}

By $\delta_i$, we denote the least in-degree  of vertices of $V(\Gamma_i)\setminus W_{i,i}$ in $\Gamma_i$. That is $$\delta_i={\rm min}\{d_{\Gamma_i}^-(v)\ |\ {v\in V(\Gamma_i)}, v\notin W_{i,i}\}.$$

For the aim to  distinguish $G_1$ from each other $G_i, i\not=1$,  we hope to prove that $\delta_1< \delta_i$ for each $i\not=1$.
 \begin{lemma} Let $\Gamma_i$ and $\delta_i$ be as defined. Then $\delta_i\geq n-2$, the equality holds if and only if one of the following cases happens: (i) $i=1$; (ii) $i=0$ and $m=5$; (iii)  $i=2$ and $m=5$.
\end{lemma}
\begin{proof}

The following claim will be used frequently.

\noindent {\sf Claim 1.}\ \ Suppose $i+k+1\leq s\leq i+2k$, $v\in W_{i,s}$. Then $|N_{\Gamma_i}^-(v)\cap W_{i,t}|\geq n-2$ for $s<t\leq i+2k$ and  $|N_{\Gamma_i}^-(v)\cap W_{i,i}|\geq \frac{n-3}{2}$.

Indeed, if $|N_{\Gamma_i}^-(v)\cap W_{i,t}|< n-2$ for $s<t\leq i+2k$, then at least two vertices in $W_{i,t}$ are out-neighbors of $v$ (noting that $W_{i,t}$ has $n-1$ vertices), which contradicts the definition of $\Gamma$. Similarly, we have $|N_{\Gamma_i}^-(v)\cap W_{i,i}|\geq \frac{n-3}{2}$.

 Now, we consider the in-degree of each vertex $v\in V(\Gamma_i)$ outside $W_{i,i}$.

  If $v$ lies in $W_{i,i+1}$, then $v=x_1$ if $i=0$ and $v=e_{i+1}$ if $i\not=0$. When $i=0$,
 $e_k\in N_{\Gamma_i}^-(x_1)$ and $ N_{\Gamma_i}^-(x_1)\supseteq  (G \setminus\{e,x\})_{k+1}$, and thus $d_{\Gamma_i}^-(v)\geq 1+(n-2)>n-2$.
 When $i\not=0$, $v=e_{i+1}$, then $N_{\Gamma_i}^-(v)\supseteq  (G\setminus\{e\})_{i+k+1}$, leading to  $d_{\Gamma_i}^-(v)>n-2$.

 If $v$ lies in $W_{i,j}$ with $i+2\leq j\leq i+k$,  then $v=e_j$.  When $j=i+k=0$, since $N_{\Gamma_i}^-(v)\supseteq (G\setminus\{e,x\})_{i+k+1}$ and $ e_{i+1}\in N_{\Gamma_i}^-(v)$, we have $d_{\Gamma_i}^-(v)>n-2$.  For all other cases,  $N_{\Gamma_i}^-(v)\supseteq  (G\setminus\{e\})_{i+k+1}$, we have  $d_{\Gamma_i}^-(v)>n-2$.

For the case when $v\in W_{i,i+k+1}$: if $k\geq 3$,  by Claim 1, $N_{\Gamma_i}^-(v)$ contains at least $n-2$ vertices from $W_{i,i+2k}$ and $W_{i,i+2k-1}$, respectively, and thus
  $d_{\Gamma_i}^-(v)\geq 2(n-2)>n-2$; if $k=2$ and $n\geq5$, since $N_{\Gamma_i}^-(v)$  contains at least $n-2$ vertices from $W_{i,i+2k}$  and it contains at least one vertex  from $W_{i,i}$ (by Claim 1),
  we have    $d_{\Gamma_i}^-(v)\geq 1+(n-2)$, as required.

  For the case when $v\in W_{i,s}$ with $i+k+2\leq s\leq i+2k-1$: $N_{\Gamma_i}^-(v)$ contains at least $n-2$ vertices from $W_{i,i+2k}$.   If $i=0$, then $v_{i+k+1} \in N_{\Gamma_i}^-(v)$, and thus  $d_{\Gamma_i}^-(v)\geq 1+(n-2)$. If $i\not=0$, then $z_{i+1}=e_{i+1}$    and thus $e_{i+1} \in N_{\Gamma_i}^-(v)$, leading to $d_{\Gamma_i}^-(v)>n-2$.

In what follows, we  suppose   $v$ is a vertex in $W_{i,s}$ with $s=i+2k$. It has been said that $W_{i,s}=(G\setminus\{w^{-1}\})_{i+2k}$, where $w=x$ if $i=1$ and $w=e$ if $i\not=1$. By Claim 1, $ N_{\Gamma_i}^-(v)$ contains at least $\frac{n-3}{2}$ elements from $W_{i,i}$. Clearly,
  $N_{\Gamma_i}^-(v)$ contains at least $\frac{n-3}{2}$ elements from $W_{i,s}$. If we can find two other in-neighbors of $v$, then we have  $d_{\Gamma_i}^-(v)>n-2$.

Firstly, we consider the case when $i=1$. In this case, $W_{i,s}=(G\setminus\{x^{-1}\})_{0}$. If $v\not=e_0$,   then  $N_{\Gamma_i}^-(v)\supseteq \{v_{i+k+1}, e_{2}\}$, proving that $d_{\Gamma_i}^-(v)>n-2$. If $v=e_0$, then   $$N_{\Gamma_i}^-(v)=(W_{1,1}\setminus\{x_1\})\bigcup (W_{1,0}\cap N_{\Gamma_i}^-(v))\bigcup\{e_{k+1}\}.$$ Since $$  |W_{1,1}\setminus\{x_1\}|=| W_{1,0}\cap N_{\Gamma_i}^-(e_0)| =\frac{n-3}{2},$$ we have   $d_{\Gamma_i}^-(e_0)=n-2$.

Now, let us consider the case when $i\not=1$.

  Suppose $m\geq 7$. If   $i+2k\not=1$, then  it follows from $N_{\Gamma_i}^-(v)\supseteq \{v_{i+2k-1}, v_{i+2k-2}\}$ that  $d_{\Gamma_i}^-(v)>n-2$.
 Assume that $i+2k=1$. If $v\not=x_{i+2k}$, then  it follows from $N_{\Gamma_i}^-(v)\supseteq \{(x^{-1}v)_{i+2k-1}, v_{i+2k-2}\}$ that  $d_{\Gamma_i}^-(v)>n-2$.
  If $v=x_{i+2k}$, noting that $z_{i+1}=e_{i+1}$, thus  $N_{\Gamma_i}^-(v)\supseteq \{e_{i+1}, x_{i+2k-2}\}$, leading to  $d_{\Gamma_i}^-(v)>n-2$.

 Next, we assume that $v\in W_{i,i+2k}$ and $m=5$.

 If $i=0$ and $v\not=x_{4}$, then  $N_{\Gamma_i}^-(v)\supseteq \{v_{3}, x_1\}$, proving that $d_{\Gamma_i}^-(v)>n-2$.  If $i=0$ and  $v=x_{4}$, then   $$N_{\Gamma_i}^-(v)=(W_{0,0}\setminus\{x_0\})\bigcup (W_{0,4}\cap N_{\Gamma_i}^-(v))\bigcup\{x_3\}.$$ Since $$|W_{0,0}\setminus\{x_0\}|=|W_{0,4}\cap N_{\Gamma_i}^-(v)|=\frac{n-3}{2},$$  we have   $d_{\Gamma_i}^-(v)=n-2$.

 If $i=2$,  and $v\not=x_{1}$, then  $N_{\Gamma_i}^-(v)\supseteq \{(x^{-1}v)_0, e_3\}$, proving that $d_{\Gamma_i}^-(v)>n-2$. If $v=x_1$, then   $$N_{\Gamma_i}^-(v)=(W_{2,2}\setminus\{x_2\})\bigcup (W_{2,1}\cap N_{\Gamma_i}^-(v))\bigcup\{e_3\}.$$ By $$|W_{2,2}\setminus\{x_2\}|=|W_{2,1}\cap N_{\Gamma_i}^-(v)|=\frac{n-3}{2},$$ we have   $d_{\Gamma_i}^-(v)=n-2$.

 If $i=3$ or $4$, then  $N_{\Gamma_i}^-(v)\supseteq \{v_{i+3}, e_{i+1}\}$, and thus   $d_{\Gamma_i}^-(v)>n-2$.
  \end{proof}
  The proof of Lemma 3.2 implies that $d_{\Gamma_i}^-(y)=n-2$, with $y\in V(\Gamma_i)\setminus W_{i,i}$,  if and only if one of the following cases happens: (1) $i=1$ and $y=e_0$; (2) $i=0, m=5$ and $y=x_4$; (3)  $i=2, m=5$ and $y=x_1$.

  If only with Lemma 3.2 in hands, one can not distinguish $\Gamma_1$ from other $\Gamma_i$ (when $m=5$, if $i=0$ or $i=2$, $\delta_i$ possibly equals to $\delta_1$). In order to prove that $\Gamma_1$ is not isomorphic to $\Gamma_0$ and $\Gamma_2$, we need further study the least out-degree of vertices in $\Gamma_i$, where $i\in\{0,1,2\}$ and $m=5$. Set $$  \pi_i={\rm min}\{d_{\Gamma_i}^+(v)\ |\  v\in V(\Gamma_i)\}.$$
   \begin{lemma} Suppose $m=5$, $n\geq 5$.   Then $\pi_0=\pi_1= \pi_2+1= \frac{ n+1}{2}$.
\end{lemma}
\begin{proof}   Suppose $i\in\{0,1,2\}$.
Since $m=5$, the vertex set of $\Gamma_i$ has a partition as:$$V(\Gamma_i)=W_{i,i}\cup W_{i,i+1} \cup W_{i,i+2}\cup W_{i,i+3}\cup W_{i,i+4}.$$ Let $v$ be a vertex in $\Gamma_i$. We consider the out-degree of $v$ in $\Gamma_i$.

If $v$ lies in $W_{i,i}$, since $N_{\Gamma_i}^+(v)$ contains at least $n-2$ vertices from $W_{i,i+3}$ and at least $n-2$ vertices from $W_{i,i+4}$, we have $d_{\Gamma_i}^+(v)\geq 2(n-2)>\frac{ n+1}{2}$.

If $v$ lies in $W_{i,i+1}$,  then $v=x_1$ if $i=0$ and $v=e_{i+1}$ if $i\in\{1,2\}$. For the case when $i=0$,
since $N_{\Gamma_i}^+(v)\supseteq  W_{0,0}$ and $N_{\Gamma_i}^+(v)\supseteq  (G\setminus\{e,x\})_4$, we have $d_{\Gamma_i}^+(v)\geq  \frac{n-1}{2}+(n-2)> \frac{ n+1}{2}$.
 When $i=1$,   $v$ exactly is $e_{2}$. In this case $N_{\Gamma_i}^+(v)\supseteq     W_{1,1}$ and $N_{\Gamma_i}^+(v)\supseteq  (G\setminus\{e,x^{-1}\})_0$,  leading to  $d_{\Gamma_i}^+(v)>\frac{ n+1}{2}$.
 When $i=2$,   $v$ is $e_{3}$. In this case $N_{\Gamma_i}^+(v)\supseteq     W_{2,2}$ and $N_{\Gamma_i}^+(v)\supseteq  \{x_1, e_4\}$,  leading to  $d_{\Gamma_i}^+(v)>\frac{ n+1}{2}$.

Suppose that  $v$ lies in $W_{i,i+2}$. In this case, $v=e_{i+2}$. When $i=0$,
 $N_{\Gamma_i}^+(v)=W_{0,0}\cup \{x_1\}$, thus   $d_{\Gamma_i}^+(v)=\frac{ n+1}{2}$.    When $i=1$,    $N_{\Gamma_i}^+(v)=W_{1,1}\cup \{e_0\}$,  proving that $d_{\Gamma_i}^+(v)=\frac{ n+1}{2}$.
 When $i=2$,   $v$ is $e_{4}$. In this case  $N_{\Gamma_i}^+(v)=W_{2,2}$,  proving that $d_{\Gamma_i}^+(v)=\frac{ n-1}{2}$.

Next, we consider the case when $v\in W_{i,j}$ with $j=i+3$ or $j=i+4$.
In this case,   $N_{\Gamma_i}^+(v)$  contains at least $ \frac{ n-3}{2}$ vertices in $W_{i,j}$. If we can find three other out-neighbors of $v$, then we have  $d_{\Gamma_i}^+(v)>\frac{ n+1}{2}$.

For the case when $v$ lies in $W_{i,i+3}=(G\setminus \{e\})_{i+3}$.  If $i=0$ and $v=x_3$, since   $N_{\Gamma_i}^+(v)\supseteq \{x_4,x_0,e_2\}$, we have  $d_{\Gamma_i}^+(v)>\frac{ n+1}{2}$. If $i=0$ and $v\not=x_3$, then  $N_{\Gamma_i}^+(v)\supseteq \{v_4,x_1,e_2\}$ it follows  $d_{\Gamma_i}^+(v)>\frac{ n+1}{2}$.
If $i=1$ and $v\not=x_4^{-1}$,  then  $N_{\Gamma_i}^+(v)\supseteq \{v_0,  e_2, e_3\}$; if $i=1$ and $v=x_4^{-1}$,  then  $N_{\Gamma_i}^+(v) =\{ e_2, e_3 \}\bigcup (N_{\Gamma_i}^+(v)\cap W_{1,4}) $, from $|N_{\Gamma_i}^+(v)\cap W_{1,4}|=\frac{n-3}{2}$, we have $d_{\Gamma_i}^+(v)=\frac{ n+1}{2}$.
If $i=2$ and $v\not=x_0^{-1}$, then    $N_{\Gamma_i}^+(v)\supseteq\{(xv)_1, e_3, e_4\} $ and thus $d_{\Gamma_i}^+(v)>\frac{ n+1}{2}$.  If $i=2$ and $v=x_0^{-1}$, then    $N_{\Gamma_i}^+(v)=\{ e_3, e_4\} \bigcup (N_{\Gamma_i}^+(v)\cap W_{i,i+3})$. Noting that $|N_{\Gamma_i}^+(x_0^{-1})\cap W_{2,0}|=\frac{n-3}{2}$, we have  $d_{\Gamma_i}^+(v)=\frac{ n+1}{2}$.

Now, only the case when  $v$ lies in $W_{i,i+4}$ is left. Suppose $v\in W_{i,i+4}$. Note that at most one element in $W_{i,i+3}$ is an in-neighbor of $v$, thus $v$ has at least $n-2$ out-neighbors in $W_{i,i+3}$.  Recalling that $v$ has at least $ \frac{ n-3}{2} $ out-neighbors   in $W_{i,i+4}$, thus we have  $d_{\Gamma_i}^+(v)>\frac{ n+1}{2}$.

The proof is completed.
\end{proof}

Suppose that  $i\in\{0,1,2\}$ and $m=5$. The proof for Lemma 3.3 implies  that $d_{\Gamma_i}^+(y)=\pi_i$, with $y\in \Gamma_i$,  if and only if one of the following cases happens: (1) $i=0$ and $y=e_2$; (2) $i=1$,  $y=e_3$ or $y=x_4^{-1}$; (3)  $i=2$ and $y=e_4$. With Lemma 3.2 and Lemma 3.3 in hands, we can distinguish $\Gamma_1$ from other $\Gamma_j$ (for the case when  $G$  admits a TRR).
 \begin{lemma}  Suppose $m\geq 7$ or $m=5$ and $n\geq 5$.  If  $\Gamma_1$ is   isomorphic to some $\Gamma_j$ with $j\in \mathbb{Z}_m$, then $j=1$ or $j=k+1$. Moreover, if   $G$  admits a TRR, then    $j=1$.
\end{lemma}
\begin{proof}  Recall that $T_{1,1}=R_1$ and $T_{l,l}=R_l$  for any $k+1\not=l\in \mathbb{Z}_m$. Thus  $\Delta_1=\Delta_l$ if $l\not=k+1$.
If $m\geq7$, since $\delta_1=n-2<\delta_j$ for any $j\not=1$,  $\Gamma_1$ is not isomorphic to $\Gamma_j$ with  $j\notin\{1,k+1\}$.   Suppose $m=5$. By comparing $\delta_1$ with $
 \delta_4$, one finds that $\Gamma_1$ is not isomorphic to   $\Gamma_4$. Since $\pi_2\not=\pi_1$, $\Gamma_1$ is not isomorphic to $\Gamma_2$. Furthermore,  $\Gamma_0$ has a unique vertex with the least out-degree,
however,  two vertices of $\Gamma_1$ (they are $e_3$ and $x_4^{-1}$) have the least out-degree, proving that $\Gamma_1$ is not isomorphic to $\Gamma_0$.

If   $G$   admits a TRR,   since $T_{1,1}=R_1$ and $T_{k+1,k+1}=R_{k+1}$, we have $\Delta_1=\Delta_{k+1}$. In this case,  from $\delta_1\not=\delta_{k+1}$ it follows that $\Gamma_1$ is not isomorphic to $\Gamma_{k+1}$. The proof   is completed.
 \end{proof}
When we proceed to prove that $\mathcal{A}_{e_1}\subseteq \mathcal{A}_{e_s}$ for any $s\not=1$ (where $\mathcal{A}$ denotes ${\rm Aut}(\Gamma)$), we   need further study the least out-degree of vertices in $V(\Gamma_i)\setminus W_{i,i}$.
For $i\in \mathbb{Z}_m$, let $\Psi_i$ be the induced subgraph of $\Gamma_i$ with vertex set $ V(\Gamma_i)\setminus W_{i,i}$. Set $   \chi _i={\rm min}\{d_{\Psi_i}^+(v)\ |\  v\in V(\Psi_i) \}.$

 \begin{lemma} Suppose $m\geq 7$ or $m=5$ and $n\geq 5$. Either $  \chi_i=0 $, or   $  \chi_i=1 $. Moreover, a vertex $v$ of $\Psi_i$ satisfies $d_{\Psi_i}^+(v)=\chi_i$ if and only if $v=e_{i+k}$.
\end{lemma}
\begin{proof} The vertex set of $\Psi_i$ has a partition as: $$V(\Psi_i)=\bigcup_{j\not=i\in\mathbb{ Z}_m} W_{i,j}.$$

For each vertex $v$ of $\Psi_i$, we consider its out-degree in $\Psi_i$.

When $v \in W_{i,i+1}$: if $i=0$, then $v=x_1$ and thus $N_{\Psi_i}^+(v)\supseteq\{x_{k+1},x^2_{2k}\}$, forcing $d_{\Psi_i}^+(v)\geq 2$;  if $i=1$, then $v=e_{2}$ and thus $N_{\Psi_i}^+(v)\supseteq\{x_{0},e_3\}$; if $i\not\in\{0,1\}$, then $v=e_{i+1}$ and thus $N_{\Psi_i}^+(v)\supseteq\{x_{i+2k},x^2_{i+2k}\}$, forcing $d_{\Psi_i}^+(v)\geq 2$.

When  $v=e_{i+j}\in W_{i,i+j}$  with $2\leq j\leq k-1$:  if $i\not=1$  it follows from $N_{\Psi_i}^+(v)\supseteq W_{i+2k}=(G\setminus\{e\})_{i+2k}$ that $d_{\Psi_i}^+(v)\geq 2$;
if $i=1$, by  $N_{\Psi_i}^+(v)\supseteq  ((G\setminus\{e,x^{-1}\})_{i+2k}\bigcup \{e_{i+k}\})$, we have $d_{\Psi_i}^+(v)\geq 2$.

When  $v=e_{i+k}\in W_{i,i+k}$: if $i+k=0$, then $N_{\Psi_i}^+(v)=\{x_1\}$ and  $d_{\Psi_i}^+(v)=1$; if $i+k=1$, then $N_{\Psi_i}^+(v)=\{e_{i+k-1}\}$ and  $d_{\Psi_i}^+(v)=1$;
if $i=1$, then $N_{\Psi_i}^+(v)=\{e_0\}$, and thus $d_{\Psi_i}^+(v)=1$; if $i=0$, then $N_{\Psi_i}^+(v)=\{x_1\}$, and thus $d_{\Psi_i}^+(v)=1$. If $i+k\notin\{0,1\}$ and $i\not\in\{0,1\}$, then $N_{\Psi_i}^+(v)=\emptyset$ and thus  $d_{\Psi_i}^+(v)=0$.

When  $v \in W_{i,i+k+1}$: if $i+k=0, v=x_1$, then  $N_{\Psi_i}^+(v)\supseteq\{e_{i+1}, x_{i+2k}\}$, and thus  $d_{\Psi_i}^+(v)\geq 2$;
if $i+k=0, v\not=x_1$, then $N_{\Psi_i}^+(v)=\{e_{i+k}, v_{i+2k}\}$, and thus $d_{\Psi_i}^+(v)\geq 2$. Suppose $i+k\not=0$. If $i\not=0$, then $N_{\Psi_i}^+(v)=\{e_{i+1}, e_{i+k}\}$, we have  $d_{\Psi_i}^+(v)\geq 2$; if $i=0$, then $N_{\Psi_i}^+(v)=\{ e_{i+k}, v_{i+2k}\}$, we have  $d_{\Psi_i}^+(v)\geq 2$.

When  $v \in W_{i,i+j}$ with $k+2\leq j\leq 2k$: it is easy to see that $N_{\Psi_i}^+(v)$ contains at least $n-2$ vertices from $W_{i,i+k+1}$. If we can find one more out-neighbor of $v$ from other $W_{i,i+j}$ with $j\not=k+1$, then we are done. If  $v\not=e_{i+j}$,  then $e_{i+k}\in N_{\Psi_i}^+(v)$   and we are done. If $v=e_{i+j}\in  W_{i,i+j}$, then $i=1$ and $j=2k$,  thus $W_{i,i+j}=W_{1,0}=(G\setminus\{x^{-1}\})_0$ and $N_{\Psi_i}^+(v)$ contains $x_0$, forcing that  $d_{\Psi_i}^+(v)\geq 2$.
\end{proof}

For $i\in \mathbb{Z}_m$, let  $\mathcal{A}_{i}$ denote the set of automorphisms of $\Gamma$ which fixes    $G_i$ point-wise. That is to say that $$\mathcal{A}_{i}=\bigcap_{g_i\in G_i} \mathcal{A}_{g_i}.$$
Lemma 3.5 is helpful for us to establish a relationship between $\mathcal{A}_{i}$ and  $\mathcal{A}_{i+k}$.
 \begin{lemma} Suppose $m\geq 7$ or $m=5$ and $n\geq 5$. We have   $\mathcal{A}_{i}\subseteq \mathcal{A}_{i+k}$ for any $i\in \mathbb{Z}_m$.
\end{lemma}
\begin{proof} Suppose $\sigma\in \mathcal{A}_{i}$. For any given vertex $v_{i+k}\in G_{i+k}$, since $\tau=\rho_{v}\cdot \sigma\cdot \rho_{v^{-1}}$ fixes $e_i$,
from $e_i^{\tau}=e_i$ it follows that $\tau$ stabilizes $\Gamma_i$. Clearly, each vertex in $W_{i,i}$ is fixed by $\tau$, thus $\Psi_i$ is stabilized by $\tau$. In view of Lemma 3.5, $e_{i+k}$ is the unique vertex in $\Psi_i$  with the least out-degree, thus $\tau$ fixes $e_{i+k}$, which implies that $\sigma$ fixes $v_{i+k}$.
 \end{proof}
With Lemma 3.4 and Lemma 3.6 in hands, we can give a proof of Theorem 1.2 for the case when $m\geq 7$ or $m=5$ and $n\geq 5$.
\begin{theorem}  Let $G$ be a finite group with odd order $n>1$. If  $m\geq 7$ or  $m=5$ and $n\geq 5$, then $G$ admits a regular TmSR.\end{theorem}
\begin{proof} The proof is partitioned into two parts according to two different cases of $G$.
\vskip 2mm
\noindent{\sf Case 1.}  $G$ admits a TRR.

Suppose $G$ has a  TRR  of the form ${\rm Cay}(G, R)$, where $R$ is a subset of $G$ with $\frac{n-1}{2}$ elements such that $R\cap R^{-1}=\emptyset,  R\cup R^{-1}=G\setminus \{e\}$. Let $x$ be a given element in $R$ and let $\Gamma$ be the graph  as defined   at the beginning part of this section.
Then $\Gamma$ is a regular tournament, the action of ${\rm Reg}(G)$ on $V(\Gamma)$ is   semiregular, and $G$ is isomorphic to ${\rm Reg}(G)$. To  prove that $G$ has a regular  TmSR, it suffices to prove that ${\rm Aut}(\Gamma)={\rm Reg}(G)$. As before, ${\rm Aut}(\Gamma)$ is abbreviated as $\mathcal{A}$.
\vskip 2mm
\noindent{\sf Claim 1.} The set $G_1$ is stabilized  by $ \mathcal{A}$.

 Let $\sigma\in \mathcal{A}$ and $u_1 \in G_1$. Suppose   $u_1^\sigma=v_l\in G_l$ with $l\in \mathbb{Z}_m$. Then $\rho_u\cdot\sigma\cdot\rho_{v^{-1}}$ sends $e_1$ to $e_l$, which implies that $\Gamma_1$ is isomorphic to $\Gamma_l$. By Lemma 3.4, $l=1$ and thus $(G_1)^\sigma\subseteq G_1$, leading to $(G_1)^\sigma=G_1$.

Since ${\rm Reg}(G)$  is transitive on $G_1$,  by Frattini argument, $\mathcal{A} =  {\rm Reg}(G) \cdot \mathcal{A}_{e_1}$. To prove that $\Gamma$ is a TmSR
over $G$, it suffices to prove $\mathcal{A}_{e_1}$ contains only the identity automorphism. Let $\tau\in  \mathcal{A}_{e_1}$. If the following claim is proved then we are done.
\vskip 2mm
\noindent{\sf Claim 2.}   $\tau$ fixes  each vertex of $\Gamma$.

  Recall the induced subgraph with vertex set $G_1$ is a TRR over $G$,  and $\tau$ stabilizes $G_1$, thus $\tau$  fixes $G_1$ point-wise.  By applying Lemma 3.6, $\tau$ fixes  $G_{k+1}$ point-wise.
  Applying  Lemma 3.6 again, we have  $\sigma$ fixes  $G_{2k+1}$ point-wise. By applying Lemma 3.6 $l$-times,  $\tau$ fixes  $G_{lk+1}$ point-wise.
  Since $k$ and $m$ are co-prime, the cyclic group $\mathbb{Z}_m$ can be generated by $k$. Thus any given $j\in \mathbb{Z}_m$ has the form $j=lk+1$ (since $j-1$ is a multiple of $k$).
  Consequently, for any given   $j\in \mathbb{Z}_m$, $\tau$ fixes  $G_{j}$ point-wise.
Hence, $\tau$ is the identity automorphism of $\Gamma$  and we are done.
\vskip 2mm
\noindent{\sf Case 2.}  $G$ is  isomorphic to $\mathbb{Z}_3^2$ or $\mathbb{Z}_3^3$.

Let $\Gamma$ be the graph  as defined at the beginning part of this section. Recall that  $T_{k+1,k+1}=S$ and  $T_{i,i}=R$ for  all other $i\in  \mathbb{Z}_m$ with $i\not=k+1$,   for a given $i\in \mathbb{Z}_m $, when $i\not=0$, we set $T_{i,i+j}=\{e\}$,  $T_{ i+j,i}=G\setminus\{e\}$ for $ j=1, \ldots, k$, and $T_{i,i+j}=G\setminus\{e\}$, $T_{ i+j,i}= \{e\}$ for $j=k+1, \ldots, 2k$; when $i=0$, set $T_{0,1}=\{x\}, T_{1,0}=G\setminus\{x^2\}$.
\vskip 2mm
\noindent{\sf Claim 3.}  For any $\sigma\in \mathcal{A}$,  $(G_1)^\sigma\subseteq (G_1\cup G_{k+1})$.

 If there is $v_1\in G_1$ such that $v_1^\sigma=u_j\in G_j$ with $j\notin\{1,k+1\}$, then $e_1^{\rho_v\sigma \rho_{u^{-1}}}=e_j$, leading to $\Gamma_1\cong \Gamma_j$, which contradicts Lemma 3.4.
\vskip 2mm
\noindent{\sf Claim 4.}   $ G_1\cup G_{k+1} $ is stabilized by any $\sigma\in \mathcal{A}$.

Firstly, we suppose $\sigma(G_1)=G_1$. Let $$\Phi=G_2\cup G_3\cup\cdots\cup G_{k+1};\ \ \ \Phi'=G_{k+2}\cup G_{k+3}\cup\cdots\cup G_{2k+1}.$$ Note that a vertex $v$ of $\Gamma$ has exactly one out-neighbor in $G_1$ if and only if $v$ is a vertex of $\Phi'$.
 Since $\sigma$ stabilizes $G_1$, we have $\sigma$ stabilizes $\Phi'$. Moreover, $\sigma$ stabilizes $\Phi$. Observe that a vertex $u$ in $\Phi$ has exactly $k$ out-neighbors in $\Phi'$ if and only if $u$ is a vertex in $G_{k+1}$.
 Consequently, $\sigma$ stabilizes $G_{k+1}$ and we are done.

  Secondly,  assume  there is some $v_1\in G_1$ such that $v_1^\sigma\in G_{k+1}$. Since $e_1^{\rho_v\sigma}\in G_{k+1}$, by replacing $\rho_v\sigma$ with $\sigma$, we may directly assume that $e_1^\sigma\in G_{k+1}$.
 Suppose $e_1^\sigma=u_{k+1}\in G_{k+1}$. We aim to prove  $(G_{k+1})^\sigma\subseteq( G_1\cup G_{k+1})$. If there is $v_{k+1}\in G_{k+1}$ such that $v_{k+1}^\sigma=w_s\in G_s$ with $s\notin\{1,k+1\}$, then by $(e_1)^{\sigma\rho_{u^{-1}}}=e_{k+1}$, $\Gamma_1\cong \Gamma_{k+1}$,
and by $(e_{k+1})^{\rho_{v}\sigma\rho_{w^{-1}}}=e_s$,   $\Gamma_{k+1}\cong\Gamma_s$, which leads to  $\Gamma_1\cong \Gamma_{s}$, which is a contradiction to Lemma 3.4. Hence, $(G_{k+1})^\sigma\subseteq( G_1\cup G_{k+1})$.
Consequently, $(G_1\cup G_{k+1})^\sigma\subseteq( G_1\cup G_{k+1})$, and we are done.

Let $\Gamma^2$ be the induced subgraph of $\Gamma$ with vertex set $G_1\cup G_{k+1}$ and let  $\sigma\in \mathcal{A}$.
By Claim 4, the restriction of $\sigma$ to $\Gamma^2$ is an automorphism of $\Gamma^2$. However, $\Gamma^2$ has ${\rm Reg}(G)$ as its automorphism group (see [4], Page
5). Thus there is $g\in G$ such that $\sigma\cdot\rho_g$ fixes $G_1\cup G_{k+1}$ point-wise. Particularly,  $\sigma\cdot\rho_g$ fixes $G_1$ point-wise.
By applying Lemma 3.6 $l$-times,  $\sigma\cdot\rho_g$ fixes $G_{lk+1}$ point-wise. Since each $j\in \mathbb{Z}_m$ can be written  as $j=lk+1$ for some $l$,  $\sigma\cdot\rho_g$ fixes each vertex of $\Gamma$.
Consequently, $\sigma=\rho_{z}$ for some $z\in G$, and thus $\mathcal{A}={\rm Reg}(G)\cong G$, as required.
\end{proof}
\section{Case II:    $m=5$ and $n=3$.}
\quad\quad When    $m=5$ and $n=3$, $G$ has a TRR  of the form ${\rm Cay}(G,R)$, where $R$ is a subset of $G$ with exactly one element that is not $e$.    Assume that  $G=\{e, x,x^2\}$ and $R=\{x\}$.  In this case, $T_{i,j}$  are   chosen in the following way for  constructing $\Gamma$:  $T_{i,i}=\{x\}$ for  each $i\in  \mathbb{Z}_5$;  if  $0\not=i\in \mathbb{Z}_5,$    we set $T_{i,i+1}=T_{i,i+2}=T_{i+3,i}=T_{i+4,i}=\{e\}$ and $T_{i+1,i}=T_{i+2,i}=T_{i, i+3}=T_{i,i+4}=\{x,x^2\}$;   particularly, we set $T_{0,1}=\{x\}, T_{1,0}= \{e,x\}$, $T_{0,2}=T_{3,0}=T_{4,0}=\{e\}$,  $T_{2,0}=T_{0,3}=T_{0,4}=\{x,x^2\}$.

It is easy to see that $\Gamma$ is a regular tournament. We aim to prove that ${\rm Aut}(\Gamma)\cong G$ and the action of ${\rm Aut}(\Gamma)$ on $V(\Gamma)$ is semi-simple with $5$ orbits.

As before, denote by $\Gamma_i$ the induced subgraph of $\Gamma$ with vertex set $N_\Gamma^+(e_i)$. Let $\Psi_i$ be the induced subgraph of $\Gamma_i$ obtained from $\Gamma_i$ by deleting the vertex $x_i$. Let  $\pi_i$ (resp., $\chi_i$) be the least out-degree of vertices in $\Gamma_i$ (in $\Psi_i$). It seems somewhat difficult to distinguish $\Gamma_1$ from other $\Gamma_i$ by using the least in-degree of vertices in $\Gamma_i$. However, we can distinguish $\Gamma_2$ from other $\Gamma_i$ by comparing the least out-degree of vertices in $\Gamma_i$.

 \begin{lemma} The indexes $\pi_i$ and $\chi_i$ are as just defined.  Then\\
 (i) \   $\pi_2=1$ and $\pi_i=2$ for $i\not=2$. \\
 (ii) \  $\chi_i=1$ if $i\not=2$ and $\chi_2=0$,  and a vertex $w$ of $\Psi_i$ satisfies $d_{\Psi_i}^+(w)=\chi_i$ if and only if $w=e_{i+2}$.
\end{lemma}
\begin{proof} Let $W_{i,j}=V(\Gamma_i)\cap G_j$. The vertex set of $\Gamma_i$ has a partition as:$$V(\Gamma_i)=W_{i,i}\cup W_{i,i+1} \cup W_{i,i+2}\cup W_{i,i+3}\cup W_{i,i+4},$$
 where   $W_{i,i}=x_i$ for  each $i\in \mathbb{Z}_5$; $W_{i,i+1}=\{e_{i+1}\}$ if $i\not=0$ and $W_{0,1}=\{x_1\}$;   $W_{i,i+2}=\{e_{i+2}\}$;
  $W_{i,i+3} =\{x_{i+3},x_{i+3}^2\}$, and $W_{i,i+4}=\{x_{i+4},x_{i+4}^{2}\}$ if   $i\not=1$,   $W_{i,i+4}=\{e_{0},x_{0}\}$ if $i=1$.

 For an arbitrary vertex $v$ of $\Gamma_i$, we study its out-degree in $\Gamma_i$. Moreover, if $v$ lies in $\Psi_i$, we also study the out-degree of $v$ in $\Psi_i$.

If $v$ lies in $W_{i,i}$, i.e., $v=x_i$, $N_{\Gamma_i}^+(v)$ contains exactly one  vertex from  each $W_{i,j}$ with $j\not=i$, and thus $d_{\Gamma_i}^+(v)=4$.

If $v$ lies in $W_{i,i+1}$,  then $v=x_1$ if $i=0$ and $v=e_{i+1}$ for other cases. When $i=0$,
since ${\Psi_i}^+(v)\supseteq\{x_3,x_4^2\}$, we have  $d_{\Gamma_i}^+(v)\geq d_{\Psi_i}^+(v)>2$; when $i=1$,  by  $N_{\Psi_i}^+(e_{2})\supseteq\{ e_3, x_0 \}$,  we have $d_{\Gamma_i}^+(v)\geq d_{\Psi_i}^+(v)>2$;  when $i\in\{2,3,4\}$,   $N_{\Psi_i}^+(v)\supseteq\{x_{i+4},x_{i+4}^2\}$, proving that   $d_{\Gamma_i}^+(v)\geq d_{\Psi_i}^+(v)>2$.

Suppose  $v$ lies in $W_{i,i+2}$. In this case, $v=e_{i+2}$. When $i=2$,   $N_{\Gamma_i}^+(v)=\{x_2\}$,  proving that $d_{\Gamma_i}^+(v)=1$ and  $d_{\Psi_i}^+(v)=0$. For all other $i$ with $i\not=2$, one can see that $d_{\Gamma_i}^+(v)=d_{\Psi_i}^+(v)+1=2$. Indeed, if $i=0$,  $N_{\Gamma_i}^+(v)=  \{x_0,x_1\}$ and $N_{\Psi_i}^+(v)=  \{x_1\}$;  if $i=1$,  then $N_{\Gamma_i}^+(v)=\{x_1, e_0\}$ and $N_{\Psi_i}^+(v)=  \{e_0\}$; if $i=3$, then $N_{\Gamma_i}^+(v)=\{x_3,x_1\}$ and $N_{\Psi_i}^+(v)=  \{x_1\}$; if $i=4$,   then  $N_{\Gamma_i}^+(v)=\{x_4,e_0\}$ and  $N_{\Psi_i}^+(v)=  \{e_0\}$.

For the case when $v$ lies in $W_{i,i+3}=\{x_{i+3},x_{i+3}^2\}$, to prove that  $d_{\Gamma_i}^+(v)\geq d_{\Psi_i}^+(v)>2$, it suffices to list  two out-neighbors of $v$ in $\Psi_i$.
\vskip 2mm
\noindent If $i=0$ and $v=x_3$, then $N_{\Psi_i}^+(v)\supseteq \{x_4,,e_2\}$;
if $i=0$ and $v=x_3^2$, then  $N_{\Psi_i}^+(v)\supseteq \{x_1,e_2\}$.\\
If $i=1$,  $v=x_4$ or $v=x_4^{2}$,  then  $N_{\Psi_i}^+(v) =\{ e_2, e_3\}$.\\
If $i=2$,  $v=x_0$ or $v=x_0^{2}$,  then    $N_{\Psi_i}^+(v)\supseteq\{e_3, e_4\} $.\\
If $i=3$ and $v=x_1$,  then  $N_{\Gamma_i}^+(v) =\{ x_2,e_4 \}$; if $i=3$ and $v=x_1^{2}$,  then  $N_{\Psi_i}^+(v)\supseteq \{ e_4, e_0\}$.\\
If $i=4$, $v=x_2$ or $v=x_2^2$,  then    $N_{\Psi_i}^+(v)=\{ e_0, e_1\}$.

When $v$ lies in $W_{i,i+4}$, to achieve the inequality $d_{\Gamma_i}^+(v)\geq d_{\Psi_i}^+(v)>2$, we only need to find two out-neighbors of $v$ in $\Psi_i$, where $W_{i,i+4}=\{x_{i+4},x_{i+4}^2\}$ if $i\not=1$, and  $W_{i,i+4}=\{e_0,x_0\}$ if $i=1$. \vskip 2mm
\noindent  If $i=0$ and $v=x_4$, then $N_{\Psi_i}^+(v)\supseteq \{x_1,e_2\}$; if $i=0$ and $v=x_4^2$, then  $N_{\Psi_i}^+(v)=\{e_2, x_3\}$.\\
If $i=1$ and $v=e_0$,  then  $N_{\Psi_i}^+(v)\supseteq \{x_4, e_2\}$; if $i=1$ and $v=x_0$,  then  $N_{\Psi_i}^+(v) =\{ x_4^2, e_3 \}$.\\
If $i=2$ and $v=x_1$ or $v=x_1^2$, then     $N_{\Psi_i}^+(v)\supseteq\{ x_0^2, e_4\} $.\\
If $i=3$ and $v=x_2$,  then  $N_{\Psi_i}^+(v) =\{ x_1^2, e_0 \}$; if $i=3$ and $v=x_2^{2}$,  then  $N_{\Psi_i}^+(v)\supseteq \{x_1, e_0\}$.\\
If $i=4$ and $v=x_3$, then    $N_{\Psi_i}^+(v)=\{x_2^2, e_1\}$; if $i=4$ and $v=x_3^2$, then    $N_{\Psi_i}^+(v)\supseteq\{ x_2, e_1\} $.

The proof is completed.
\end{proof}
\begin{theorem}  Let $G$ be a finite group with $3$ elements. If    $m=5$, then $G$ admits a regular TmSR.\end{theorem}
\begin{proof} Let $\Gamma$ be the regular tournament  defined as above.
  Since the action of ${\rm Reg}(G)$ on $V(\Gamma)$ is   semi-regular with $5$ orbits and $G$ is isomorphic to ${\rm Reg}(G)$, we only need to prove that ${\rm Aut}(\Gamma)={\rm Reg}(G)$. As before, ${\rm Aut}(\Gamma)$ is abbreviated as $\mathcal{A}$.

  Let $\tau\in \mathcal{A}$ and $u_2\in G_2$. Suppose   $u_2^\tau=v_l\in G_l$ with $l\in \mathbb{Z}_5$. Then $\rho_u\cdot\tau\cdot\rho_{v^{-1}}$ sends $e_2$ to $e_l$, which implies that $\Gamma_2$ is isomorphic to $\Gamma_l$, and thus $\pi_2=\pi_l$. By Lemma 4.1, $l=2$ and thus $(G_2)^\tau\subseteq G_2$, leading to $(G_2)^\tau=G_2$.

Since ${\rm Reg}(G)$  is transitive on $G_2$,  by Frattini's argument, $\mathcal{A} =  {\rm Reg}(G) \cdot \mathcal{A}_{e_2}$. To prove that $\Gamma$ is a TmSR
over $G$, it suffices to prove $\mathcal{A}_{e_2}$ contains only the identity automorphism.

Let $\sigma\in  \mathcal{A}_{e_2}$. We hope to prove that $\sigma$ fixes each vertex of $\Gamma$.
   Since the induced subgraph with vertex set $G_2$ is a TRR over $G$,   $\sigma$  fixes $G_2$ point-wise.

        For any given vertex $v_{4}\in G_{4}$, since $\sigma'= \rho_{v}\cdot \sigma\cdot \rho_{v^{-1}}$ fixes $e_2$, it stabilizes $\Gamma_2$. Since $\sigma$ fixes $x_2v_2$, $\sigma'$  fixes $x_2$. Consequently, $\sigma'$ stabilizes $\Psi_2$.  However, $e_{4}$ is the unique vertex in $\Psi_2$  with the least out-degree, $e_{4}$  must be fixed by   $\sigma'$, which leads to $v_{4}^\sigma=v_{4}$. Hence, $\sigma$ fixes $G_4$ point-wise.

    Now, we consider the action of $\sigma$ on $G_1$.  For any given vertex $v_{1}\in G_{1}$, since $\sigma'= \rho_{v}\cdot \sigma\cdot \rho_{v^{-1}}$ fixes $e_4$ and $x_4$, it stabilizes $\Psi_4$. However, $e_{1}$ is the unique vertex in $\Psi_4$  with the least out-degree, $e_{1}$  is fixed by   $\sigma'$, which leads to $v_{1}^\sigma=v_{1}$. Consequently, $\sigma$ fixes $G_1$ point-wise.

     We  continue the action of $\sigma$ on $G_3$.  For any given vertex $v_{3}\in G_{3}$, since $\sigma'= \rho_{v}\cdot \sigma\cdot \rho_{v^{-1}}$ fixes $e_1$ and $x_1$, it stabilizes $\Psi_1$. However, $e_{3}$ is the unique vertex in $\Psi_1$  with the least out-degree, $e_{3}$  is fixed by   $\sigma'$, which leads to $v_{3}^\sigma=v_{3}$. Consequently, $\sigma$ fixes $G_3$ point-wise.

     Finally, we study the action of $\sigma$ on $G_0$.  For a given vertex $v_{0}\in G_{0}$, since $\sigma'= \rho_{v}\cdot \sigma\cdot \rho_{v^{-1}}$ fixes $e_3$ and $x_3$, it stabilizes $\Psi_3$. However, $e_{0}$ is the unique vertex in $\Psi_3$  with the least out-degree, $e_{0}$  is fixed by   $\sigma'$, which leads to $v_{0}^\sigma=v_{0}$. Consequently, $\sigma$ fixes $G_0$ point-wise.

     Therefore, $\sigma$ is the identity automorphism of $\Gamma$ and thus  ${\rm Aut}(\Gamma)={\rm Reg}(G)$.
     \end{proof}

\section{Case III:   $m=3$.}
\quad\quad If $m=3$,  $T_{i,j}$  are   chosen in the following way for  constructing $\Gamma$.
Set  $T_{0,0}=T_{1,1}=R^{-1}$ if   $G$ has a TRR, or $T_{0,0}=T_{1,1}=S$  if $G \cong \mathbb{Z}_3^2$ or $G   \cong  \mathbb{Z}_3^3$,  $T_{2,2}=R$,  $T_{0,1}=\{x\}, T_{1,0}=G\setminus\{x^{-1}\}$, $T_{1,2}=T_{2,0}=\{e\}$, $T_{2,1}=T_{0,2}=G\setminus\{e\}$, where $R, S$ and $x$ are given in Proposition 2.1, 2.2 and 2.3 (if $G$ has a TRR, $x$ is a given element in $R$).

Our goal is to prove that ${\rm Aut}(\Gamma)$ is isomorphic to $G$ and the action of ${\rm Aut}(\Gamma)$ on $V(\Gamma)$ is semi-regular with exactly $3$-orbits.

Let $\Gamma_i$ be the induced subgraph of $\Gamma$ with vertex set $N_{\Gamma_i}^+(e_i)$. The vertex set of $\Gamma_i$ has a partition as $$V(\Gamma_i)=W_{i,i}\cup W_{i,i+1}\cup W_{i,i+2},$$
where $W_{i,j}=N_{\Gamma}^+(e_i)\cap G_j$ for $i\leq j\leq i+2$. Clearly, for $i\in \{0,1\},$ $W_{i,i}=R_i^{-1}$ if $G$ has a TRR and $W_{i,i}=S_i$ if $G \cong \mathbb{Z}_3^2$ or $G\cong  \mathbb{Z}_3^3$, $W_{2,2}=R_2$; $W_{0,1}=\{x_1\}$, $W_{0,2}=(G\setminus\{e\})_2$;
$W_{1,2}=\{e_2\}$, $W_{1,0}=(G\setminus x^{-1})_0$; $W_{2,0}=\{e_0\}$, $W_{2,1}=(G\setminus\{e\})_1$.

 A key step is to prove  $\Gamma_1$ is not isomorphic to $\Gamma_0$ or $\Gamma_2$. To achieve the goal, we  study the least out-degree of  vertices in $\Gamma_i$.
As in Section 3, let $$  \pi_i={\rm min}\{d_{\Gamma_i}^+(v)\ |\  v\in V(\Gamma_i)\}.$$
\begin{lemma}  $\pi_i\geq \frac{n-3}{2}$, with equality if and only if $i=1$.  Moreover, a vertex $v$ of $\Gamma_1$ satisfies  $d_{\Gamma_i}^+(v)=\frac{n-3}{2}$ if and only if $v=e_0$.
\end{lemma}
\begin{proof}
Let $v$ be a vertex of $\Gamma_i$.
If $v\in W_{i,i}$, since at most one vertex in $W_{i,i+2}$ is an in-neighbor of $v$, we have $d_{\Gamma_i}^+(v)\geq n-2>\frac{n-3}{2}$.

If $v\in W_{i,i+1}$, noting that all vertices in $W_{i,i}$ are out-neighbors of $v$, we have $d_{\Gamma_i}^+(v)\geq  \frac{n-1}{2}$.

Next, suppose $v\in W_{i,i+2}$.

 If $i=0$ and $v=x_2$, note that $(e_2, x_2)$ is an arc of $\Gamma$ and $e_2\notin W_{0,2}$, $x_2$ has   $\frac{n-1}{2}$ out-neighbors in $W_{0,2}$, thus we have $d_{\Gamma_i}^+(v)\geq  \frac{n-1}{2}$.
If $i=0$ and $v\not=x_2$, note that  $(v, x_1)$ is an arc of $\Gamma_0$, $v$ has at least  $\frac{n-3}{2}$ out-neighbors in $W_{0,2}$, thus we have $d_{\Gamma_i}^+(v)> \frac{n-3}{2}$.
Consequently, $\pi_0> \frac{n-3}{2}$.

 If $i=2$ and $v=x_1$, note that $(x_1, x_2)$ is an arc of $\Gamma_2$ and $x_1$ has  at least $\frac{n-3}{2}$ out-neighbors in $W_{2,1}$, thus we have $d_{\Gamma_i}^+(v)> \frac{n-3}{2}$.
If $i=2$ and $v\not=x_1$, note that  $(v, e_0)$ is an arc of $\Gamma_2$, $v$ has at least  $\frac{n-3}{2}$ out-neighbors in $W_{2,1}$, thus we have $d_{\Gamma_i}^+(v)> \frac{n-3}{2}$.
Consequently, $\pi_2> \frac{n-3}{2}$.

 If $i=1$ and $v\not=e_0$, then  $(v, e_2)$ is an arc of $\Gamma_1$, thus  $d_{\Gamma_i}^+(v)> \frac{n-3}{2}$ (noting that  $v$ has at least  $\frac{n-3}{2}$ out-neighbors in $W_{1,0}$).
 If $i=1$ and $v=e_0$, note that $(e_0, x_0^{-1})$ is an arc of $\Gamma$ and $x_0^{-1}$ is not in $W_{1,0}$, thus   $e_0$ has  exactly $\frac{n-3}{2}$ out-neighbors in $W_{1,0}$.
It is easy to see that $e_0$ has no out-neighbors in $W_{1,1}$. Consequently,  $d_{\Gamma_1}^+(v)= \frac{n-3}{2}$ and thus   $\pi_1=\frac{n-3}{2}$.

Since $e_0$ is the unique vertex in $\Gamma_1$ with out-degree $\frac{n-3}{2}$, the second assertion follows.
  \end{proof}

 \begin{theorem}  Let $G$ be a finite group with odd order. If    $m=3$, then $G$ admits a regular TmSR.\end{theorem}
\begin{proof} Let $\Gamma$ be the regular tournament  defined as above.
 We aim to prove that ${\rm Aut}(\Gamma)={\rm Reg}(G)$. As before, ${\rm Aut}(\Gamma)$ is abbreviated as $\mathcal{A}$.

By Lemma 5.1, it is clear that $G_1$ is stabilized by $\mathcal{A}$. To see this, let $\tau\in \mathcal{A}$ and $v_1\in G_1$. Suppose   $v_1^\tau=w_l\in G_l$ with $l\in \mathbb{Z}_3$. Then $\rho_v\cdot\tau\cdot\rho_{w^{-1}}$ sends $e_1$ to $e_l$, which implies that $\Gamma_1$ is isomorphic to $\Gamma_l$, and thus $\pi_1=\pi_l$. By Lemma 5.1, $l=1$ and thus $(G_1)^\tau\subseteq G_1$, leading to $(G_1)^\tau=G_1$.

\vskip 2mm
\noindent{\sf Case 1.} $ G$ admits a TRR.

Since ${\rm Reg}(G)$  is transitive on $G_1$,  by Frattini argument, $\mathcal{A} =  {\rm Reg}(G) \cdot \mathcal{A}_{e_1}$. To prove that $\Gamma$ is a regular TmSR
over $G$, it suffices to prove $\mathcal{A}_{e_1}$ contains only the identity automorphism.
By Lemma 5.1, it is easy to see that $\mathcal{A}_{e_1}\subseteq \mathcal{A}_{e_0}$. Indeed, if    $\tau\in  \mathcal{A}_{e_1}$, then $\tau$ stabilizes $\Gamma_1$ and thus $\tau$ fixes $e_0$ (since $e_0$ is the unique vertex in $\Gamma_1$ with the least out-degree).

Let $\sigma\in \mathcal{A}_{e_1}$.  We hope to prove $\sigma$ fixes each vertex of $\Gamma$.
  Recall $G_1$ is stabilized by $\sigma$ and the induced subgraph with vertex set $G_1$ is a TRR over $G$, thus $\sigma$  fixes $G_1$ point-wise.

 For an arbitrary vertex $v_0\in G_0$, from $ \rho_{v}\cdot\sigma\cdot\rho_{v^{-1}} \in\mathcal{A}_{e_1} $ it follows that $ \rho_{v}\cdot\sigma\cdot\rho_{v^{-1}}\in\mathcal{A}_{e_0}$, which implies  that $(v_0)^{\sigma}=v_0$.
 Furthermore, $\sigma$ stabilizes $G_2$.  For an arbitrary vertex $w_2\in G_2$, since $w_2$ is the unique out-neighbor of $w_1\in G_1$, it follows from $(w_1)^{\sigma}=w_1$ that $(w_2)^{\sigma}=w_2$.
 Consequently, $\sigma$ is the identity automorphism and we are done.

\vskip 2mm
\noindent{\sf Case 2.}  $G$ is  isomorphic to $\mathbb{Z}_3^2$ or $\mathbb{Z}_3^2$.

  For any $\tau\in \mathcal{A}$,  since $(G_1)^\tau=G_1$, it is easy to see that $G_2\cup G_0$ is stabilized by $\tau$.

Let $\Gamma^2$ be the induced subgraph of $\Gamma$ with vertex set $G_2\cup G_{0}$ and let  $\sigma\in \mathcal{A}$.
The restriction of $\sigma$ to $\Gamma^2$ is an automorphism of $\Gamma_2$. However, $\Gamma^2$ has ${\rm Reg}(G)$ as its automorphism group (by Proposition 2.2, 2.3). Thus there is some $g\in G$ such that $\sigma\cdot\rho_g$ fixes $G_2\cup G_{0}$ point-wise. For an arbitrary vertex $w_1\in G_1$, since $w_1$ is the unique out-neighbor of $(x^{-1}w)_0$ and  $(x^{-1}w)_0$ is fixed by  $\sigma\cdot\rho_g$, it follows that $w_1$ is fixed by  $\sigma\cdot\rho_g$.
Consequently,  $\sigma\cdot\rho_g$ is the identity automorphism of $\Gamma$ and thus $\sigma=\rho_{z}$ for some $z\in G$, which implies that  $\mathcal{A}={\rm Reg}(G)\cong G$, as required.
 \end{proof}
\section{Conclusion}
\quad\quad In this article, we prove that a finite group of odd order $n>1$ has a regular TmSR if $m\geq 3$ is odd. Thus all nonidentity finite groups of odd order with a regular
TmSR are classified.  We remind that: the case when $G=\mathbb{Z}_1$ has only the identity element is left. The purpose of  ``graph representation of finite groups" is: View an element of an abstract finite group $G$ as an automorphism of a graph. Since an automorphism of a graph seems  more concrete than an element of a group, ``graph representation" is helpful for  understanding   the elements of an abstract group.
In view of this point, the problem  ``Whether  $\mathbb{Z}_1$ admits a  TmSR" seems  missing some importance since the element in $\mathbb{Z}_1$ is clear.
 Even so, we feel a little regret  for leaving out the trivial group $\mathbb{Z}_1$.
 Du \cite{du} has observed that   $\mathbb{Z}_1$  has no regular T3SR.
It is not difficult to see that    $\mathbb{Z}_1$  has no regular T5SR. To see this, we refer the reader to the following lemma.

\begin{lemma} A regular  tournament with $5$ vertices has nontrivial automorphism.
\end{lemma}
\begin{proof} Let $\Gamma$ be a  regular tournament with $5$ vertices. Thus each vertex of $\Gamma$ must have $2$ out-neighbors.
Suppose $v_1$ is a given vertex of $\Gamma$. Assume, without loss of generality, $v_2, v_3$ are  out-neighbors of $v_1$ and $v_4, v_5$ are  in-neighbors of $v_1$.
Suppose $(v_2, v_3)$ is an arc of $\Gamma$ (if $(v_3, v_2)$ is an arc, we relabel $v_3$ as $v_2$ and relabel $v_2$ as $v_3$). Since $v_3$ already have two in-neighbors  $v_1, v_2$, the rest two vertices $v_4, v_5$ must be
out-neighbors of $v_3$.

If $(v_2, v_4)$ is an arc, to assure the regularity of $\Gamma$, $(v_4,v_5)$ and $(v_5,v_2)$ must be arcs of $\Gamma$. In this case, $\Gamma$ is the Cayley graph ${\rm Cay}(\mathbb{Z}_5: X)$ with  $X=\{1,2\}$. Clearly, $\Gamma$ has an  automorphism  which sends $v_i$ to $v_{i+1}$, for $1\leq i\leq 4$, and which sends $v_5$ to $v_1$.

If $(v_4, v_2)$ is an arc, to assure the regularity of $\Gamma$, $(v_5,v_4)$ and $(v_2,v_5)$ must be arcs of $\Gamma$. In this case, $\Gamma$  has an  automorphism  of the form $(v_1 \ v_2\ v_3\ v_5\ v_4)$ (where the automorphism is written as a permutation on $V(\Gamma)$).
\end{proof}

If $m\geq 7$ is odd, we believe that   $\mathbb{Z}_1$    has no regular TmSR. However, when $m$ is large enough, the proof is somewhat elusive.

\end{document}